\documentclass[12pt]{book}

\usepackage{amsmath}
\usepackage{amssymb}
\usepackage{array}
\usepackage{delarray}
\usepackage{amsfonts}
\usepackage{amsthm}
\usepackage[english]{babel}
\usepackage{enumerate}
\usepackage{fancyhdr}
\usepackage{slashbox}
\usepackage{graphicx}
\usepackage{float}

\RequirePackage[OT1]{fontenc}

\theoremstyle{plain}
\newtheorem{theorem}{Theorem}

\newtheorem{lem}{Lemma}

\newcommand{\Bea} {\begin{eqnarray*}}
\newcommand{\Eea} {\end{eqnarray*}}
\newcommand{\R}{\mathbb{R}}
\newcommand{\N}{\mathbb{N}}

\newcommand{\W}{\mathbb{W}}
\newcommand{\E}{\mathbb{E}}
\newcommand{\p}{\mathbb{P}}

\begin{document}

\renewcommand{\baselinestretch}{1.2}
\markright{
}

\markboth{\hfill{\footnotesize\rm Thanh Mai Pham Ngoc}\hfill}
{\hfill {\footnotesize\rm Statistical minimax approach of the
Hausdorff moment problem} \hfill}
\renewcommand{\thefootnote}{}
$\ $\par \fontsize{10.95}{14pt plus.8pt minus .6pt}\selectfont
\vspace{0.8pc} \centerline{\large\bf Statistical minimax approach of
the} \vspace{2pt} \centerline{\large\bf Hausdorff moment problem}
\vspace{.4cm} \centerline{Thanh Mai Pham Ngoc} \vspace{.4cm}
\centerline{\it Universit{\'e} Pierre et Marie Curie, Paris VI}
\vspace{.55cm} \fontsize{9}{11.5pt plus.8pt minus .6pt}\selectfont

\begin{quotation}

\noindent \textit{Abstract}: The purpose of this paper is to study
the problem of estimating a compactly supported function from noisy
observations of its moments in an ill-posed inverse problem
framework. We provide a statistical approach to the famous Hausdorff
classical moment problem. We prove a lower bound on the rate of
convergence of the mean integrated squared error and provide an
estimator which attains minimax rate over the corresponding
smoothness classes.
\par
\vspace{9pt} \noindent {\it Key words and phrases:} {Minimax
estimation}, Statistical inverse problems, {Problem of the moments},
{Legendre polynomials}, {Kullback-Leibler information}
\par \vspace{9pt} \noindent {\it MSC 2000 Subject Classification}
62G07 44A60  45Q05
\par
\end{quotation}\par

\fontsize{10.95}{14pt plus.8pt minus .6pt}\selectfont

\setcounter{chapter}{1}
\setcounter{equation}{0} 
\noindent {\bf 1. Introduction}

The classical moment problem can be stated as follows: it consists
in getting some information about a distribution $\mu$ from the
knowledge of its moments $\int x^{k}d\mu(x)$. This problem has been
largely investigated in many mathematical topics, among others, in
operator theory, mathematical physics, inverse spectral theory,
probability theory, inverse problems, numerical analysis but also in
a wide variety of settings in physics and engineering such as
quantum spin models, speech processing, radio astronomy (see for
instance Kay and Marple (1981), Lang and McClellan (1983) $\dots$).
We may cite the classical and pioneer books in the field (see
Akhiezer (1965), Shohat and Tamarkin (1943)) which put emphasis on
the existence aspect of the solution and its uniqueness. According
to the support of the distribution of interest, one may refer to one
of the three types of classical moment problems: the Hamburger
moment problem whose support of $\mu$ is the whole real line, the
Stieljes problem on $[0,+\infty)$ and finally the Hausdorff problem
on a bounded interval. In this paper, we shall focus on the last
issue under the inverse problem angle.

The Hausdorff moment problem which dates back to 1921 (see Hausdorff
(1921)) occupies a central place in the field of inverse problems
and has been an object of great interest in literature since then.
For instance, the particular case when only a finite number of
moments are known which is called the Truncated Hausdorff moment
problem, has recently aroused much attention (see Talenti (1987),
Fasino and Inglese (1996), Inglese (1995), Tagliani (2002)). Another
interest and aspect of the Hausdorff moment problem lies in its very
close link to the inversion of the Laplace transform when this
latter is given at equidistant points ( see for instance Brianzi and
Frontini (1991), Dung, Huy, Quan and Trong (2006), Tagliani (2001),
Tagliani and Velazquez (2003)). In fact, by a simple change of
variable, the problem of Laplace transform inversion is equivalent
to the Hausdorff moment problem. More recently, Ang, Gorenflo, Le
and Trong (2002) have presented the Hausdorff moment problem under
the angle of ill posed problems, in a sense that solutions do not
depend continuously on the given data. Nonetheless, until now, as
far as we know, the statistical approach which consists in assuming
that the noise is stochastic has been very little put forward and
rarely raised.

We consider in this paper a statistical point of view of the
Hausdorff moment problem. We aim at reconstructing an unknown
function from noisy measurements of its moments on a symmetric
bounded interval in a statistical inverse problem framework. In
reality, it is barely impossible to measure moments with having any
corruption. Without loss of generality we may and will suppose that
$[-a,a]=[-1,1]$. In practise, observations of moments appear
especially in quantum physics (see for instance Ash and Mc Donald
(2003) and Mead and Papanicolaou (1983)), image analysis (see Teague
(1980)), engineering mechanics (see Athanassoulis and Gavriliadis
(2002)) and all the references therein). We shall give now some
concrete examples of the application of the problem of moments, in
which either one can have access directly to moments to reconstruct
the unknown function at stake or
the moment problem appears in the line of reasoning.  \\
Examples. 1. In quantum physics, a big issue aims at reconstructing
a positive density of states in harmonic solids from its moments on
a finite interval. It is explained how to measure these moments (p
2410 and
2411 in Mead and Papanicolaou (1983)) and (p 3078 and 3079 in Gaspard and Cyrot-Lackmann (1973)).\\
2. In the field of quantum gravity, the transition probabilities for
a Markov chain related to the causal set approach to modeling
discrete theories of quantum gravity satisfy a moment problem (see
Ash and McDonald (2003)). One has to measure those probabilities. \\
3. In the context of non classical moment problem, we may cite the
following example (see Ang, Nhan and Thanh (1999)). It deals with
the determination of the shape of an object in the interior of the
Earth by gravimetric methods. The density of that object differs
from the density of the surrounding medium. Assuming a flat earth
model, the problem consists in finding a curve $x \rightarrow
\sigma(x)$ in the half plane $0\leq \sigma(x)<H, 0\leq x\leq 1$,
$\sigma(x)$ satisfying a non linear integral equation of the first
kind of the form $$
\frac{1}{2\pi}\int_{0}^{1}\frac{H-\sigma(\xi)}{(x-\xi)^2+(H-\sigma(\xi))^2}d\xi=f(x),$$
 where $f(x)$ is a given function. The nonlinear integral equation can be approximated by the following
 linear integral equation in $\varphi$:
 $$ \int_{0}^{1} \frac{\varphi(\xi)}{(M+x+\xi)^2}d\xi=2\pi f(-M-x),$$
  with $\varphi(x)=H-\sigma(x)$, $M$ is large enough and $x\geq 0$. By taking $x=1,2,\dots,n,\dots$, we get the following equivalent
  moment problem:
  $$ \int_{0}^{1} \frac{\varphi(\xi)}{(M+n+\xi)^2}d\xi=\mu_{n},$$
  where $\mu_{n}=2\pi f(-M-n), \, n=1,2,\dots$
  \\

In the first examples cited above, the authors aim at estimating an
unknown density from measurements of its moments. The following
paper presents the case of the reconstruction of an unknown function
which could be in particular a density of probability, when one has
its moments corrupted with some white noise. This approach
constitutes one among other ill-posed inverse problems point of
view.
 \\

Estimation in statistics using moments has
already been put forward (see Mnatsakanov and Ryumgaart (2005)), but
the approach there was based on empirical moments and empirical
processes, the results were expressed in terms of weak convergence
whereas our paper is built upon an ill-posed inverse approach with
minimax results.
\\  The estimation procedure we use is based on the
expansion of the unknown function through the basis of Legendre
polynomials and an orthogonal series method. We establish an upper
bound and a lower bound on the estimation accuracy of the procedure
showing that it is optimal in a minimax sense. We show that the
achieved rate is only of logarithmic order. This fact has already
been underlined by Goldenshluger and Spokoiny (2004). In their
paper, Goldenshluger and Spokoiny (2004) tackled the problem of
reconstructing a planar convex set from noisy geometric moments
observations. They pointed out that in view of reconstructing a
planar region from noisy measurements of moments, the upper bound
was only in the order of logarithmic rate. The lower bound has not
been proved. In a second part, they consider reconstruction from
Legendre moments to get faster rates of convergence. Legendre
moments can be observed in the context of shapes reconstruction. In
our present work, instead of considering a planar region, we deal
with functions belonging to a Sobolev scale and we stay focused on
the classic moments with respect to the monomials $x^{k}$. Moreover,
recently in the context of long-memory processes obtained by
aggregation of independant parameter AR(1) processes and in view of
estimating the density of the underlying random parameter, Leipus,
Oppenheim, Philippe and Viano (2006) had to deal with a problem of
moments. They obtained very slow logarithmic rate but without
showing that this could be the best possible. In a certain way, our
minimax results provide a piece of answer.

One might question this chronic slow rate which seems inherent to
moment problems. In fact, the underlying problem lies in the non
orthogonal nature of the monomials $x^{k}$. They actually hamper the
convergence rate to be improved for bringing a small amount of
information. This remark is highlighted in our proof of the upper
bound.

This paper is organized as follows: in section 2 we introduce the
model and the estimator of the unknown function and we finally state
the two theorems. Section 3 contains the proofs. The last section is
an appendix in which we prove some useful
inequalities about binomial coefficients.\\

\par
\setcounter{chapter}{1}
\setcounter{equation}{0} 
\noindent{\bf 2. Statement of the problem}\\

\noindent \textbf{ 2.1 The model.} First of all, let us recall an
usual statistical
framework of ill-posed inverse problems (see Mathe and Pereverzev (2001)):\\
 Let $A: H
\longrightarrow H$ be a known linear operator on a Hilbert space
$H$. The
problem is to estimate an unknown function $f \in H$ from indirect observations \\
\begin{equation}\label{modele_de_base}
Y=Af+\varepsilon\xi,
\end{equation}
where $\varepsilon$ is the amplitude of the noise. It is supposed to
be a small positive parameter which tends to $0$, $\xi$ is assumed
to be a zero-mean Gaussian random process indexed by $H$ on a
probability space
$(\Omega,\mathcal{A},\bold P) $.\\
In the Hausdorff moment problem, the operator $A_{k}$ which
determines the way the observations are indirect is defined by:
\begin{equation}\label{operateur}
A_{k}(f)=\int_{-1}^{1}x^{k}f(x)dx,\quad  k\in \N.
\end{equation}
Let us now state the Hausdorff moment problem model. From the
equation (\ref{modele_de_base}) and (\ref{operateur}), we derive the
following sequence of moments observations pertubated by a
stochastic noise:
\begin{equation}\label{modele} y_ {k}=\mu_{k}+\varepsilon \xi_{k} \quad
k=0,1,\dots.
\end{equation}
 where $\varepsilon$ is the noise level, it is supposed to tend to $0$, $\xi_{k}$ are assumed to be i.i.d standard
Gaussian random variables and $\mu_{k}$ are the moments of the
unknown function $f$ given by:
$$\mu_{k}=\int_{-1}^{1}x^{k}f(x)dx \quad k=0,1,\dots.$$
The assumption that the noise is modelled as i.i.d standard Gaussian
random variables is a natural modeling and has already been used
(see Rodriguez and Seatzu (1993), Brianzi and Frontini (1991),
Goldenshluger and Spokoiny (2004)). It is a standard assumption in
statistical inverse problems. Moreover, we may consider
(\ref{modele}) as a first approach
of the problem at hand. \\
The objective is to estimate the unknown function $f$ supported on
the bounded interval $[-1,1]$ from noisy observations of its moments
in the model (\ref{modele}) which can be assimilated to a gaussian
white noise model but using a non orthonormal basis and which
constitutes an inverse problem setting. The solution of this problem
is unique in the case $k=\infty$ (see Bertero, De Mol and Pike
(1985)), so that the statistical problem of recovering $f$ from
(\ref{modele}) is relevant. In the particular case of estimating a
density of probability, the probability measure having
density $f$ is unique (see Feller (1968, Chap.7)). \\
The use of the Legendre polynomials in the Hausdorff classical
moment problem in order to approximate the unknown measure is
 a completely standard and quite natural procedure (see Ang, Gorenflo, Le and
Trong (2002), Bertero, De Mol and Pike (1985), Papoulis (1956),
Teague (1980), Goldenshluger and Spokoiny (2004)) as those
polynomials directly result from the Gram-Schmidt orthonormalization
of the family $\{x^{k}\},\, k=0, 1, \dots$ Consequently, expanding
the function $f$ to be estimated in the basis of Legendre
polynomials
falls naturally and fits the problem's nature. \\
Denote $\beta_{n,j}$ the coefficients of the
normalized Legendre polynomial of
 degree $n$:
 $$P_{n}(x)=\sum_{j=0}^{n}\beta_{n,j}\,x^{j}$$
 By considering the
Legendre polynomials, if we multiply both sides of the model
(\ref{modele}) by the coefficients $\beta_{k,j}$ we get the
following model (for the proof see Lemma \ref{equivalence_modeles}.
in Appendix):
\begin{equation}\label{modele_hetero}
\tilde y_{k}=\theta_{k}+\varepsilon\sigma_{k}\xi_{k}
\end{equation}
where $\sigma_{k}^{2}=\sum_{j=0}^{k} \beta_{k,j}^{2}$,\, $\tilde
y_{k}=\sum_{j=0}^{k} \beta_{k,j}y_{j}$,\,
$\theta_{k}=\sum_{j=0}^{k}\beta_{k,j}\mu_{j}=\int_{-1}^{1}f(x)P_{k}(x)dx$
and  $\xi_{k}$ are i.i.d standard Gaussian random variables. So the
model (\ref{modele}) is equivalent to the model
(\ref{modele_hetero}). The model (\ref{modele_hetero}) will be used
for the proof of the lower bound.\\Before going any further, we can
make a remark at this stage concerning the model
(\ref{modele_hetero}) which is an heteroscedastic gaussian sequence
space model. Depending on the asymptotic behavior of the intensity
noise $\sigma_{n}^{2}$ one may characterize the nature of the
problem's ill-posedness (see Cavalier, Golubev, Lepski and Tsybakov
(2004)). Here, in our case, $\sigma_{n}^{2} \geq \frac{1}{4}4^{n}$
(see Lemma \ref{lemmebino2}. from Appendix) and hence tends to
infinity exponentially. We may say that we are dealing with a
severely ill-posed problem with log-rates
.\\
\newline We assume that $f$ belongs to
the Sobolev space $W^{r}_{2}$ defined by:
$$W^{r}_{2}=\{f\in L^{2}[-1,1]:
\sum_{k} k^{2r}|\theta_{k}|^{2}<\infty\}$$ where
$\theta_{k}=\int_{-1}^{1}f(x)P_{k}(x) dx $ is the Legendre Fourier
coefficient and $P_{k}$ denotes the normalized Legendre polynomial
of degree $k$. Note that we consider more general functions than
densities of probability which are included in this smoothness class.\\
 Sobolev spaces associated with various kinds of underlying orthonormal
basis constitute quite standard smoothness assumption classes in
classical ill-posed problems (see for instance Mair and Ruymgaart
(1996), Mathe and Pereverzev (2002), Goldenshluger and Pereverzev
(2000)). In the Hausdorff moment problem, the underlying basis is
the Legendre
polynomials.\\
Let us now give some highlights of the Sobolev space $W^{r}_{2}$
regarding Legendre polynomials. Rafal'son (1968) and Tomin (1973) have shown in the more
general case of Jacobi polynomials (and thus in the particular case
of Legendre polynomials we are considering here) that Sobolev space
$W^{r}_{2}$ consists of all functions $f$ which have their
derivatives $f',f'',\dots,f^{(r-1)}$ being
absolutely continuous on each interval $[a,b] \subset (-1,1)$ (see also Mathe and Pereverzev (2002)).\\
\newline \noindent \textbf{ 2.2 The estimation procedure.}
Let us define now the estimator of $f$. This latter is induced by an
orthogonal series method through the Legendre polynomials.\\
Any function in $L^{2}[-1,1]$ has an expansion:\\
 $$ f(x)=\sum_{k=0}^{\infty}\theta _{k}P_{k}(x) \quad \textrm{with}
 \quad \theta_{k}=\int_{-1}^{1}f(x)P_{k}(x) dx. $$
 The problem of estimating $f$ reduces to estimation of the
 sequence $\{\theta_{k}\}_{k=1}^{+\infty}$ for Legendre polynomials
 form a complete orthogonal function system in $L^{2}[-1,1]$.

 We have
 $$\theta_{k}=\int_{-1}^{1}f(x)P_{k}(x)dx=\sum_{j=0}^{k}\beta_{k,j}\int_{-1}^{1}f(x)x^{j}dx=\sum_{j=0}^{k}\beta_{k,j}\mu_{j}.$$
 This leads us to consider the following estimator of $\theta$:
 $$\hat{\theta}_{k}=\sum_{j=0}^{k}\beta_{k,j}y_{j}$$
 and hence the estimator $\hat{f}_{N}$ of $f$:
 $$\hat{f}_{N}(x)=\sum_{k=0}^{N} \hat{\theta}_{k}P_{k}(x)=\sum_{k=0}^{N}\sum_{j=0}^{k}\beta_{k,j}y_{j}P_{k}(x)$$
 where $y_{j}$ is given by (\ref{modele}) and $N$ is an integer to be
 properly selected later.

 The mean integrated square error of the estimator $\hat{f}_{N}$ is:
 $$\E_{f}\|\hat{f}_{N}-f\|^{2},$$
 where $\E_{f}$ denotes the expectation w.r.t the distribution
 of the data in the model (\ref{modele}) and for a function $g \in
 L^{2}[-1,1]$,
 $$ \|g\|=\bigg (\int_{-1}^{1}g^{2}(x)dx \bigg ) ^{1/2}.$$
  In this paper we shall consider the problem of estimating $f$ using the mean
 integrated square risk in the model (\ref{modele}).

 We state now the two results of the paper. The first theorem establishes an upper bound.
\begin{theorem}\label{theo1}
For $\alpha >0$, define the integer $N=\lfloor \alpha
\log{(1/\varepsilon)}\rfloor$. Then we have
$$ \sup_{f \in \W^{r}_{2}} \E_{f} \|\hat{f}_{N}-f\|^{2}\leq
C[\log(1/\varepsilon)]^{-2r},$$  where $C$ is an absolute positive
constant and $\lfloor \cdot \rfloor$ denotes the floor function.
\end{theorem}
\noindent We recall that the floor function is defined by: $\lfloor
x \rfloor=\max\{n \in
\mathbb Z | \, n \leq x \}$. \\
 \noindent The second theorem provides a lower bound.
\begin{theorem}\label{theo2}
We have
$$ \inf_{\hat{f}}\sup_{f \in \W^{r}_{2}} \E_{f}\|
\hat{f}-f\|^{2}\geq c[\log(1/\varepsilon)]^{-2r},$$  where $c$ is a
positive constant which depends only on $r$ and the infimum is taken
over all estimators $\hat{f}$.
\end{theorem}

\par

\setcounter{chapter}{3}
\setcounter{equation}{0} 
\noindent {\bf 3. Proofs}\\

\noindent \textbf{ 3.1 Proof of Theorem \ref{theo1}.} For the
following proof, we consider the genuine model (\ref{modele}). By
the usual MISE decomposition which involves a variance term and a
bias term, we get
$$\E_{f}\|\hat{f}_{N}-f\|^{2}=\E_{f}\sum_{k=0}^{N}(\hat{\theta}_{k}-\theta_{k})^{2}+\sum_{k\geq N+1}\theta_{k}^{2}$$
but \Bea \E_{f}\sum_{k=0}^{N}(\hat{\theta}_{k}-\theta_{k})^{2} &=& \E_{f}\sum_{k=0}^{N}(\sum_{j=0}^{k}\beta_{k,j}(y_{j}-\mu_{j}))^{2}\\
&=& \varepsilon^{2}
\E_{f}\sum_{k=0}^{N}(\sum_{j=0}^{k}\beta_{k,j}\xi_{j})^{2}
\Eea \\
and since $\xi_{j} \stackrel {iid}{\sim} N(0,1)$, it follows that
\Bea \E_{f}\|\hat{f}_{N}-f\|^{2}&=&\varepsilon^{2}
\sum_{k=0}^{N}\sum_{j=0}^{k}\beta_{k,j}^{2} + \sum_{k\geq N+1}
\theta_{k}^{2} \\
&=& V_{N} + B_{N}^{2} \Eea \\
 We first deal with the variance term $V_{N}$. To this end,
we have to upper bound the sum of the squared coefficients of the
normalized Legendre polynomial of degree $k$. Set
$\sigma_{k}^{2}=\sum_{j=0}^{k}\beta_{k,j}^{2}$. An explicit form of
$P_{k}(x)$ is given by (see Abramowitz and Stegun (1970)):
$$P_{k}(x)=\bigg(\frac{2k+1}{2}\bigg)^{1/2}\frac{1}{2^{k}}\sum_{j=0}^{[k/2]}(-1)^{k}\binom{k}{j}\binom{2k-2j}{k}x^{k-2j},$$
where $[\cdot]$ denotes the integer part and $\binom{k}{j}$ denotes
the binomial coefficient, \\$\binom{k}{j}=\frac{k!}{(k-j)!j!}$. This
involves
 \begin{eqnarray} \sigma_{k}^{2} &=&
\frac{2k+1}{2}\frac{1}{4^{k}}\sum_{j=0}^{[k/2]}\bigg\{{\binom{k}{j}\binom{2k-2j}{k}\bigg\}}^{2}
\nonumber \\ \label{ineq_sigma1}&\leq&
\frac{2k+1}{2}\frac{1}{4^{k}}{\bigg\{\binom{2k}{k}\bigg\}}^{2}\sum_{j=0}^{[k/2]}\bigg\{{\binom{k}{j}\bigg\}}^{2}
\\ \label{ineq_sigma2} &\leq&
\frac{2k+1}{2}\frac{1}{4^{k}}{\bigg\{\binom{2k}{k}\bigg\}}^{2}(2^{k})^{2},
\end{eqnarray}\\
the inequality (\ref{ineq_sigma1}) is due to the fact that for
$0\leq j \leq [k/2] $, we have $\binom{2k-2j}{k}\leq \binom{2k}{k}$.
As for (\ref{ineq_sigma2}), we have
$\sum_{j=0}^{[k/2]}\big\{{\binom{k}{j}\big\}}^{2} \leq
\big\{\sum_{j=0}^{k}{\binom{k}{j}\big\}}^{2}$ and it is well known
that $\big\{\sum_{j=0}^{k}{\binom{k}{j}\big\}}^{2}=(2^k)^2$.
\\By using now that ${\{\binom{2k}{k}\}}^{2}\leq \frac{4^{2k}}{\sqrt{k}}$ (see
Lemma \ref{lemmebino1}. from Appendix) we have
$$\sigma_{k}^{2} \leq \frac{2k+1}{2}\frac{4^{2k}}{\sqrt{k}}$$
which yields
$$ V_{N} \leq C{\varepsilon}^{2}N^{3/2}4^{2N},$$
where $C>0$ denotes an absolute positive constant.

Now, it remains to upper bound the bias term $B_{N}^{2}$.
 \Bea
 B_{N}^{2}&=&\sum_{k\geq
N+1}\theta_{k}^{2} \\&=&\sum_{k\geq
N+1}\frac{k^{2r}}{k^{2r}}\theta_{k}^{2} \\&\leq& N^{-2r}\sum_{k=1}^
{\infty}{k^{2r}}\theta_{k}^{2} \Eea Since the function $f$ belongs
to the space $W^{r}_{2}$, $\sum_{k} k^{2r}|\theta_{k}|^{2}<\infty$,
we get
$$B_{N}^{2}=\mathcal{O}(N^{-2r})$$
Finally we have the upper bound for the MISE:
\begin{equation}\label{balance}
\E_{f}\|\hat{f}_{N}-f\|^{2} \leq
C{\varepsilon}^{2}N^{3/2}4^{2N}+C'N^{-2r}
\end{equation}
 At last, it remains to choose the optimal $N$ which will minimize the expression (\ref{balance}). This $N$ is obtained by
 equalizing the upper bounds of the bias and the variance term, namely:
 $$C{\varepsilon}^{2}N^{3/2}4^{2N}=C'N^{-2r},$$
 as $4^{2N}\gg N^{2r+3/2}$, consequently
$ N\asymp\log(\frac{1}{\varepsilon^2})$. Once one plugs
$N\asymp\log(\frac{1}{\varepsilon^2})$ in
(\ref{balance}), the desired result of the Theorem 1. follows. \\

\par

\noindent \textbf{3.2 Proof of Theorem 2.} From now on, to prove the
lower bound and for practical reasons, we shall consider the model
(\ref{modele_hetero}) which constitutes a heteroscedastic gaussian
sequence space model. We recall that the equivalence between the
models
(\ref{modele}) and (\ref{modele_hetero}) is proved in Lemma \ref{equivalence_modeles}. in Appendix.\\
A successful approach and standard tool to obtain lower bounds for
minimax risk consists in specifying a subproblem namely constructing
a subset of functions based on the observations
(\ref{modele_hetero}). Then we lean on the application of the
following particular version of Fano's lemma (see Birg\'e and
Massart (2001)) which will allow us to evaluate the difficulty of
the specified subproblem and will give us a lower bound for the MISE
associated to this subproblem.
\\ One crucial point in the Fano's lemma is the use of the
Kullback-Leibler divergence $K(\p_{1},\p_{0})$ between two
probability distributions $\p_{1}$ and $\p_{0}$ defined by:
\begin{displaymath}
K(\p_{1}, \p_{0})=\left\{
\begin{array}{ll} \int_{\R} \log(\frac{p_{1}(x)}{p_{0}(x)})p_{1}(x)dx & \textrm{if} \; \p_{1}\ll \p_{0} \\
+\infty & \textrm{otherwise}
\end{array}\right.
\end{displaymath}
Here's the version of Fano's lemma, we are going to exploit:
\begin{lem}
Let $\eta$ be a strictly positive real number, $\mathcal C$ be a
finite set of elements $\{f_{0},\dots,f_{M}\}$ on $\R$ with
$|\mathcal C|\geq
6$ and $\{P_{j}\}_{j\in \mathcal{C}}$ a set of probability measures indexed by $\mathcal{C}$ such that : \\
(i)\, $\|f_{i}-f_{j}\| \geq \eta >0,\quad \forall \, 0\leq i<j \leq M$.\\
(ii) $\p_{j} \ll \p_{0}, \quad \forall \, j=1,\dots,M,$ and
$$K(\p_{j},\p_{0})\leq H < \log M$$
then for any estimator $\hat{f}$ and any nondecreasing function
$\ell$
$$ \sup_{f\in \mathcal C} \E_{f} \big[ \ell( \|\hat{f}-f\|) \big ]
\geq \ell(\frac{\eta}{2}) \bigg [1-\bigg (\frac{2}{3} \vee
\frac{H}{\log M}\bigg ) \bigg ].$$
\end{lem}

$\newline$
 First of all, we have to construct an appropriate set of functions $\mathcal E$.
  We are going to define $\mathcal E$ as a set of functions of the following type
\begin{eqnarray}\mathcal E=\bigg \{f_{\delta}\in W_{2}^{r}:
f_{\delta}=\mathbf{1}_{[-1,1]}(\frac{c_{0}}{m^{(4r+3)/2}}\sum_{k=m}^{2m-1}\delta_{k}
k^{(2r+2)/2} P_{k}), \; \nonumber \\ \delta=(\delta_{m}, \dots,
\delta_{2m-1}) \in \Delta=\{0,1\}^{m}\bigg \}\nonumber
\end{eqnarray}

\noindent We verify that $f_{\delta}$ belongs to $W_{2}^{r}$. In
this aim, we have to calculate the Legendre-Fourier coefficients
associated with the function $f_{\delta}$:
\begin{eqnarray}
\theta_{\delta \, l}&=&\int_{-1}^{1}f_{\delta}(x)P_{l}(x)dx \nonumber \\
&=& \label{coefflegendre} \left \{ \begin{array}{ll}
\frac{c_{0}}{m^{(4r+3)/2}}\cdot l^{(2r+2)/2}\cdot \delta_{l}
& \textrm{if}\; l\in[m,2m-1] \\
0 & \textrm{else} \end{array} \right.
\end{eqnarray}
hence
\begin{eqnarray}
\sum_{k=0}^{+\infty}k^{2r}\theta_{\delta \, k}^{2}&=&
\frac{c_{0}^{2}}{m^{4r+3}}\sum_{k=m}^{2m-1}k^{2r}k^{2r+2}
\delta_{k}^{2} \nonumber \\ &\leq&
\frac{c_{0}^{2}}{m^{4r+3}}\sum_{k=m}^{2m-1}k^{4r+2}\delta_{k}^{2} \nonumber \\
&\leq&
\frac{c_{0}^{2}(2m)^{4r+2}}{m^{4r+3}}\sum_{k=m}^{2m-1}\delta_{k}^{2}\leq
c_{0}^{2}2^{4r+2} < \infty, \nonumber \end{eqnarray} since
$\delta_{k}\in\{0,1\}$.

We set $\delta^{(0)}=(0,\dots,0)$ and $f_{\delta^{(0)}}\equiv
f_{0}$. The Legendre-Fourier coefficients of $f_{0}$ are null:
\begin{eqnarray}\label{coeffnull}
 \theta_{0\,l}&=& 0 \quad \forall \, l  \,
\in \N.
\end{eqnarray}

\noindent  We are now going to exhibit the suitable subset of
functions $\mathcal C$ of the Lemma 1. To this purpose, we only take
into consideration a subset of $M+1$ functions of $\mathcal E$:
$$ \mathcal C= \{f_{\delta^{(0)}},\dots,f_{\delta^{(M)}}\}$$
where $\{\delta^{(1)},\dots,\delta^{(M)}\}$ is a subset of
$\{0,1\}^{m}$. \\
We precise that $\p_{{\delta}}$ is the law of the vector of
observations $\tilde{Y}=(\tilde{y}_{1},\dots,\tilde{y}_{\infty})$ in
the model (\ref{modele_hetero}) for $f=f_{\delta}$, $\delta \in
\mathcal{C}$.

We are now going to apply Lemma 1. We first check the condition (i),
accordingly, we have to assess the distance
$\|f_{\delta^{(i)}}-f_{\delta^{(j)}}\|^{2}$. By the orthogonality of
the system $\{P_{k}\}_{k}$ and thanks to Parseval equality we get,
for $0 \leq i<j\leq M$,
\begin{eqnarray}
\|f_{\delta^{(i)}}-f_{\delta^{(j)}}\|^{2} &=&
\frac{c_{0}^{2}}{m^{4r+3}} \sum_{k=m}^{2m-1}
k^{2r+2}(\delta_{k}^{(i)}-\delta_{k}^{(j)})^{2} \nonumber \\
&\geq& \frac{c_{0}^{2}}{m^{4r+3}}\cdot m^{2r+2}\sum_{k=m}^{2m-1}
(\delta_{k}^{(i)}-\delta_{k}^{(j)})^{2} \nonumber\\
&\geq& \frac{c_{0}^{2}}{m^{2r+1}}
\sum_{k=m}^{2m-1}(\delta_{k}^{(i)}-\delta_{k}^{(j)})^{2} \nonumber \\
&=& \frac{c_{0}^{2}}{m^{2r+1}} \rho(\delta^{(i)},\delta^{(j)}),
\nonumber
\end{eqnarray}
where $\rho(\cdot,\cdot)$ is the \textit{Hamming distance}.
 We are going to resort to the Varshamov-Gilbert bound which is
 stated in the following lemma to find a lower bound of the quantity $\rho(\delta^{(i)},\delta^{(j)})$:
\begin{lem}{\rm{(Varshamov-Gilbert bound, 1962)}.}
Fix $m \geq 8$. Then there exists a subset
$\{\delta^{(0)},\dots,\delta^{(M)}\}$ of $\Delta$ such that $M\geq
2^{m/8}$ and
$$\rho(\delta^{(j)},\delta^{(k)}) \geq \frac{m}{8}, \quad \forall \, 0\leq
j <k \leq M .$$ Moreover we can always take
$\delta^{(0)}=(0,\dots,0)$.
\end{lem}
\noindent For a proof of this lemma see for instance Tsybakov
(2004), p 89.
\\Hence
$$\|f_{\delta^{(i)}}-f_{\delta^{(j)}}\|^{2} \geq
(c_{0}^{2})/(8m^{2r})\equiv \eta^2$$ \\
We are now going to check the condition (ii) in Lemma 1. and
evaluate the Kullback-Leibler divergence. It is well known (see for
instance Birg\'e and Massart (2001) p 62) that for the
Kullback-Leibler divergence in the case of a gaussian sequence space
model we have
\begin{eqnarray}\label{kullbackhetero}
K(\p_{{\delta}},\p_{{0}})&=&\frac{1}{\varepsilon^{2}}\sum_{l=1}^{\infty}\frac{|\theta_{\delta
\, l}-\theta_{0 \, l}|^{2}}{\sigma_{l}^{2}}.
\end{eqnarray}

\noindent Hence, by virtue of (\ref{coefflegendre}),
(\ref{coeffnull}) and (\ref{kullbackhetero}), the Kullback-Leibler
divergence between the two probability measures $\p_{{0}}$ and
$\p_{{\delta}}$ of the observations in the model
(\ref{modele_hetero}) associated respectively with functions
$f=f_{0}$ anf $f=f_{\delta}$ for all $\delta \in \mathcal C$
satisfies
\begin{eqnarray}\nonumber
K(\p_{{\delta}},\p_{{0}})&=&\frac{1}{\varepsilon^{2}}\frac{c_{0}^{2}}{m^{4r+3}}\sum_{l=m}^{2m-1}\frac{l^{2r+2}\delta_{l}^{2}}{\sigma_{l}^{2}}\\
&\leq&
\frac{c_{0}^{2}2^{2r+2}}{\varepsilon^{2}}\frac{m^{2r+2}}{m^{4r+3}}\sum_{l=m}^{2m-1}\frac{\delta_{l}^{2}}{\sigma_{l}^{2}}
\nonumber
\end{eqnarray}
but thanks to Lemma \ref{lemmebino2}. (see Appendix) we have
$$ \frac{1}{\sigma^{2}_{l}} \leq \frac{1}{4^{l-1}}$$
which implies
\begin{eqnarray} \nonumber
K(\p_{{\delta}},\p_{{0}})&\leq& \frac{c_{0}^{2}
2^{2r+4}}{\varepsilon^{2}}\frac{1}{m^{2r+1}4^{m}} \sum_{l=m}^{2m-1}
\delta_{l}^{2} \leq \frac{c_{0}^{2}
2^{2r+4}}{\varepsilon^{2}}\frac{1}{m^{2r+1}4^{m}}\cdot m\leq
\frac{c_{0}^{2}2^{2r+4} m}{\varepsilon^{2}4^{m}}
\end{eqnarray}
One chooses $m=\frac{1}{\log 4}\log(\frac{1}{\varepsilon^{2}})$ so
that $1/4^m=\varepsilon^2$, hence:
$$K(\p_{{\delta}},\p_{{0}}) \leq c_{0}^{2}2^{2r+4}m $$
and since $m\leq 8 \log M/\log2$ (see Lemma 2)
$$K(\p_{{\delta}},\p{_{0}}) \leq \frac{c_{0}^{2}2^{2r+7}}{\log 2}\log M$$
Eventually one can choose $c_{0}$ small enough to have
$c\equiv\frac{c_{0}^{2}2^{2r+7}}{\log2}< 1$.

Since now all the conditions of the Lemma 1. are fulfilled, we are
in position to apply its result with the loss function
$\ell(x)=x^{2}$ and $\eta=(c_{0})/(2\sqrt{2}m^{r})$. Therefore, we
derive that whatever the estimator $\hat{f}$,
$$ \sup_{f\in \mathcal C} \E_{f}[\|\hat{f}-f\|^{2}]
\geq \frac{c_{0}^{2}}{32m^{2r}}\bigg [1-\bigg (\frac{2}{3} \vee
\frac{c}{\log M}\bigg ) \bigg ]\geq \frac{c_{0}^{2}}{96m^{2r}}.$$
but from above we had $m=\frac{1}{\log
4}\log(\frac{1}{\varepsilon^{2}})$ which gives the desired lower
bound.

\hfill \\

\setcounter{chapter}{4}
\setcounter{equation}{0} 
\noindent {\bf 4. Conclusion}

The two theorems of this paper show that in the problem of
estimating a function on a compact interval from noisy moments
observations, the best rate of convergence one can achieve,
supposing Sobolev scale smoothness and considering the mean
integrated squared error, is only of logarithmic order. In a future
work, one could try to generalise this result for $L^p$ loss and may
obtain faster rate of convergence if one assumes a more restricted
smoothness class involving super smooth functions. Besides, one may
consider an heteroscedastic gaussian noise instead of a white noise
model.

\hfill \\

\setcounter{chapter}{5}
\setcounter{equation}{0} 
\noindent {\bf 5. Appendix}

\begin{lem}\label{equivalence_modeles}
The models (\ref{modele}) and (\ref{modele_hetero}) are equivalent.
\end{lem}
\par
\noindent \textbf{Proof}. We recall the model (\ref{modele}):
$$y_{j}=\mu_{j}+\varepsilon \xi_{j}=\int_{-1}^{1}f(x)x^jdx +\varepsilon \xi_{j}.$$
We are going now to multiply both sides of (\ref{modele}) by the
coefficient $\beta_{jk}$ of the Legendre polynomial:
\begin{eqnarray}
(\ref{modele})&\Leftrightarrow&
\beta_{kj}y_{j}=\beta_{kj}\int_{-1}^{1}f(x)x^j dx+ \varepsilon
\beta_{kj} \xi_{j} \nonumber \\&\Leftrightarrow&
\sum_{j=0}^{k}\beta_{kj}y_{j}=\sum_{j=0}^{k}\beta_{kj}\int_{-1}^{1}f(x)x^j
dx + \sum_{j=0}^{k}\varepsilon \beta_{kj} \xi_{j} \nonumber \\
&\Leftrightarrow&
\sum_{j=0}^{k}\beta_{kj}y_{j}=\int_{-1}^{1}f(x)\sum_{j=0}^{k}\beta_{kj}x^{j}dx+\varepsilon\sum_{j=0}^{k}
\beta_{kj}\xi_{j} \nonumber  \\&\Leftrightarrow&
\sum_{j=0}^{k}\beta_{kj}y_{j}=\int_{-1}^{1}f(x)P_{k}dx +
\varepsilon\sum_{j=0}^{k} \beta_{kj}\xi_{j} \nonumber
\end{eqnarray}
Let us set $\tilde{\xi}_{k}=\sum_{j=0}^{k} \beta_{kj}\xi_{j}$. Since
$\xi_{j}$ are i.i.d standard Gaussian random variables, the random
variable $\tilde{\xi}_{k}$ follows a normal law with zero mean and
variance equal to $\sum_{j=0}^{k} \beta_{kj}^{2}$. Hence
(\ref{modele}) is equivalent to:
\begin{equation}\nonumber
\tilde y_{k}=\theta_{k}+\varepsilon\sigma_{k}\xi_{k}
\end{equation}
where $\sigma_{k}^{2}=\sum_{j=0}^{k} \beta_{k,j}^{2}$,\, $\tilde
y_{k}=\sum_{j=0}^{k} \beta_{k,j}y_{j}$,\,
$\theta_{k}=\sum_{j=0}^{k}\beta_{k,j}\mu_{j}=\int_{-1}^{1}f(x)P_{k}(x)dx$
and  $\xi_{k}$ are i.i.d standard Gaussian random variables.

\begin{lem}\label{lemmebino1}
For all $n\geq 1$ we have:
\begin{equation}\label{inegalitesigma1}
\binom{2n}{n} \leq \frac{4^{n}}{n^{1/4}}
\end{equation}
\end{lem}
\par
\noindent \textbf{Proof}.  Let us prove (\ref{inegalitesigma1}) by
recursion on $n$. The inequality is
clearly true for $n=1$.\\
Suppose (\ref{inegalitesigma1}) true for a certain $n \geq 1$.
\begin{eqnarray} \nonumber
\binom{2(n+1)}{n+1}= \binom{2n}{n}\frac{2(2n+1)}{n+1} \leq
\frac{4^{n}}{n^{1/4}}\frac{2(2n+1)}{n+1},
\end{eqnarray}
by recursion hypothesis. It remains to prove that
\begin{equation}\label{inegalitesigma3}
 \frac{4^{n}}{n^{1/4}}\frac{2(2n+1)}{n+1} \leq
 \frac{4^{n+1}}{{(n+1)}^{1/4}}.
\end{equation}
\begin{eqnarray}
 (\ref{inegalitesigma3}) & \Longleftrightarrow&
\frac{2(2n+1)}{n^{1/4}(n+1)} \leq \frac{4}{(n+1)^{1/4}} \nonumber \\
& \Longleftrightarrow & \frac{n+1}{n}\big (\frac{2n+1}{n+1}
\big)^{4} \leq 2^{4} \nonumber \\
& \Longleftrightarrow & (n+\frac{1}{2})^{1/4} \leq n(n+1)^{3},
\nonumber
\end{eqnarray}
which is true because we have $(n+\frac{1}{2})^{1/4} \leq
(n+\frac{1}{2})^{3}(n+1)$ and $(n+\frac{1}{2})^{3}\leq n(n+1)^{2}$
since $\frac{1}{8} \leq n^{2}/2+n/4$. This completes the proof.

\begin{lem}\label{lemmebino2}
For all $n \geq 1$, we have:
\begin{equation}\label{inegalitesigma2}\sigma_{n}^{2} \geq
4^{n-1}\end{equation}
\end{lem}
where $\sigma_{n}$ is defined in (\ref{modele_hetero}).
\\
\par
\noindent \textbf{Proof.} Firstly, let us recall the value of the
noise intensity $ \sigma_{n}^{2}$:
\begin{eqnarray}
\sigma_{n}^{2}&=&\frac{2n+1}{2}\frac{1}{4^{n}}\sum_{j=0}^{[n/2]}\bigg\{{\binom{n}{j}\binom{2n-2j}{n}\bigg\}}^{2} \nonumber \\
&\geq& \frac{n}{4^{n}}{\binom{2n}{n}}^{2}. \nonumber
\end{eqnarray}
And so, in order to prove (\ref{inegalitesigma2}) it remains to
prove that
$$ \binom{2n}{n} \geq \frac{4^{n}}{2\sqrt{n}} \qquad n \geq 1.$$
We again use a recursion on $n$.\\
The inequality (\ref{inegalitesigma2}) is clear for $n=1$. We
suppose the property true for a certain $n \geq 1$ and we shall
prove it at the rank $(n+1)$.
\begin{eqnarray}
\binom{2(n+1)}{n+1}&=&\binom{2n}{n}\frac{2(2n+1)}{n+1}\nonumber \\
&\geq& \frac{4^{n}}{2\sqrt{n}}\frac{2(2n+1)}{n+1}\nonumber \\
\label{lastinequality} &>& \frac{4^{n+1}}{2\sqrt{n+1}}
\end{eqnarray}
the inequality (\ref{lastinequality}) is true because it is
equivalent to $4n^{2}+4n+1>4n^{2}+4n$ what we always have.
\\
\par

\noindent{\large\bf Acknowledgements}

The author wishes to thank the referees for their valuable
criticisms and suggestions and Dominique Picard for interesting
discussions, leading to the present improved version of the paper.
\\
\par

\noindent{\large\bf References}
\begin{description}

\item
Abramowitz, M. and Stegun, I. A. (1970). \textit{Handbook of
Mathematical Functions with Formulas, Graphs, and Mathematical
tables}. Dover Publications, New York.

\item
Akhiezer, N. I. (1965). \textit{The classical Moment Problem and
Some Related Questions in Analysis}. Oliver and Boyd, Edinburgh.

\item
Ang, D. D., Gorenflo, R., Le, V. K. and Trong, D. D. (2002).
Moment Theory and Some inverse Problems in Potential Theory and Heat
Conduction. \textit{Lectures Notes in Mathematics} \textbf{1792},
Springer, Berlin.

\item
Ang, D. D., Nhan, V. N.; Thanh, D. N. (1999). A nonlinear integral
equation of gravimetry: Uniqueness and approximation by linear
moments. \textit{Vietnam J. Math.}, \textbf{27}(1), 61-67. (1999).

\item
Athanassoulis G.A. and Gavriliadis P. N. (2002). The truncated
Hausdorff moment problem solved by using kernel density functions.
\textit{Probabilistic Engineering Mechanics}, \textbf{17}, 273-291.

\item
Ash, A. and Mc Donald P. (2003). Moments problem and the causal set
approach to quantum gravity. \textit{Journal of Mathematical
Physics}, \textbf{44}(4), 1666-1678.

\item
Birg\'e, L. and Massart, P. (2001). Gaussian model selection.
\textit{J. Eur. Math. Soc}, \textbf{3}, 203-268.

\item
Bertero, M., De Mol, C., and Pike, E. R. (1985). Linear inverse
problems with discrete data. I: General formulation and singular
system analysis. \textit{Inverse Problems}, \textbf{1}, 301-330.

\item
Brianzi, P. and Frontini, M. (1991). On the regulized inversion of
the Laplace transform. \textit{Inverse Problems}, \textbf{7},
355-368.

\item
Cavalier, L., Golubev, Y., Lepski, O., and Tsybakov, A. B.
(2004). Block thresholding and sharp adaptative estimation in
severely ill-posed inverse problems. \textit{Theory Probab. Appl.},
\textbf{48}(3), 426-446.

\item
Dung, N., Huy, N. V., Quan, P. H., Trong, D. D. (2006). A
Hausdorff-like moment problem and the inversion of the Laplace
transform. \textit{Mathematische Nachrichten}, \textbf{279}(11),
1147-1158.

\item
Fasino, D. and Inglese, G. (1996). Recovering a probability
density from a finite number of moments and local a priori
information. \textit{Rend. Instituto. Mat. Univ. Trieste},
\textbf{28}(1), 184-200.

\item
Feller, W. (1968). \textit{An introduction to Probability Theory and
its Applications. Vol II. 2nd ed.} Wiley Series in Probability and
Mathematical Statistics, New York.

\item
Gaspard, J.P. and Cyrot-Lackman, F. (1973). Density of states from
moments. Application to the impurity band. \textit{J. Phys. C: Solid
State Phys}, \textbf{6}, 3077-3096.

\item
Goldenshluger, A. and Pereverzev, S. (2000). Adaptive estimation of
linear functionals in Hilbert scales from indirect white noise
observations. \textit{Probab. Theory Relat. Fields}, \textbf{118},
169-186.

\item
Goldenshluger, A. and Spokoiny, V. (2004). On the
shape-from-moments problem and recovering edges from noisy Radon
data. \textit{Probab. Theory Relat. Fields}, \textbf{128}, 123-140.

\item
Hausdorff, F. (1921). Summationsmethoden und Momentfolgen, I.
\textit{Mathematische Zeitschrift}, \textbf{9}, 74-109.

\item
Inglese, G. (1995). Christoffel functions and
finite moments problems. \textit{Inverse problems}, \textbf{11},
949-960.

\item
Kay, M. and Marple, S. L. (1981). Spectrum analysis-a modern
perspective. \textit{ Proc. IEEE}, \textbf{69}, 1380-1419.

\item
Lang, S. W. and McClellan, J. H. (1983). Spectral estimation for
sensor arrays. \textit{IEEE Trans. Acoust. Speech Signal Process},
\textbf{31}, 349-358.

\item
Leipus, R., Oppenheim, G., Philippe, A., Viano, M. C. (2006).
Orthogonal series density estimation in a disaggregation scheme.
\textit{J. Stat. Plan. Inf.}, \textbf{138}(8), 2547-2571.

\item
Mair B., Ruymgaart, F. H. (1996). Statistical estimation in Hilbert
scales. \textit{SIAM. J. Appl. Math.}, \textbf{56}, 1424-1444.

\item
Mathe, P. and Pereverzev, S. V. (2001). Optimal discretization of
inverse problems in Hilbert scales. Regularization and
self-regularization of projections methods. \textit{SIAM J. Numer.
Anal.} \textbf{38}(6), 1999-2021.

\item
Mathe, P. and Pereverzev, S. V. (2002). Stable summation of
orthogonal series with noisy coefficients. \textit{J. Approx.
Theory} \textbf{ 118}(1), 66-80.

\item
Mead, L. R. and Papanicolaou, N. (1984). Maximum entropy in the
problem of moments. \textit{J. Math. Phys.} \textbf{25}(8),
2404-2417.

\item
Papoulis, A. (1956). A new method of inversion of the Laplace
transform. \textit{Q. Appl. Math.}, \textbf{14}, 405-414.

\item
Rafal'son S.Z. (1968). The approximation of functions by
Fourier-Jacobi sums.  \textit{Izv. Vyss. Ucebn. Zaved. Matematika },
\textbf{4}, 54-62.

\item
Rodriguez, G. and Seatzu S. (1992). On the solution of the finite
moment problem. \textit{J. Math. Anal. Appl.}, \textbf{171} 321-333.

\item
Rodriguez, G. and Seatzu S. (1993). Approximation methods for finite
moment problem. \textit{Numerical algorithms}, \textbf{5}(8)
391-405.

\item
Shohat, J. A. and Tamarkin, J. D. (1943). \textit{The problem
of moments}. American Mathematical Society, New York.

\item
Tagliani, A. (2001). Numerical inversion of Laplace transform on the
real line of probability density functions. \textit{Applied
mathematics and computation}, \textbf{123}(3), 285-299.

\item
Tagliani, A. (2002). Entropy estimate of probability densities
having assigned moments: Hausdorff case. \textit{Appl. Math. Lett.},
\textbf{15}(3), 309-314.

\item
Tagliani, A. and Velazquez Y. (2003). Numerical inversion of the
Laplace transform via fractional moments. \textit{Applied
mathematics and computation}, \textbf{143}(1),

\item
Talenti, G. (1987). Recovering a function from a finite number of
moments. \textit{Inverse Problems}, \textbf{3}, 501-517.

\item
Teague, M. R. (1980). Image analysis via the general theory of
moments. \textit{J. Opt. Soc. Am.}, \textbf{70}, 920-924.

\item
Tomin, N. G (1973). An application of the interpolation of linear
operators to questions of the convergence of series of Fourier
coefficients with respect to classical orthogonal polynomials.
\textit{Dokl. Akad. Nauk. SSSR}, \textbf{212}, 1074-1077.

\item
Tsybakov, A. B. (2004). \textit{Introduction \`a
l'estimation non-param\'etrique}. Springer-Verlag. Berlin.
\end{description}

\noindent Laboratoire de Probabilit\'es et Mod\`eles al\'eatoires,
UMR 7599, Universit\'e Paris 6, case 188, 4, Pl. Jussieu, F-75252
Paris Cedex 5, France. \\
E-mail: thanh.pham\_ngoc@upmc.fr

\end{document}